\documentclass[12pt,leqno]{article}
\usepackage{amsmath, amsthm, amsfonts, amssymb, color}
\usepackage{mathrsfs}
\theoremstyle{definition}

\newcommand{\scr}[1]{\mathscr #1}
\definecolor{wco}{rgb}{0.5,0.2,0.3}

\numberwithin{equation}{section} \theoremstyle{remark}

\newcommand{\ua}{\uparrow}

\title{{\bf Derivative Formula and Gradient Estimates for  Gruschin Type Semigroups}\footnote{Supported in
 part by NNSFC(11131003), SRFDP, the Laboratory of Mathematical and  Complex Systems and the Fundamental Research Funds for the Central Universities.}
}
\author{
{\bf Feng-Yu Wang}\\
\footnotesize{School of Mathematical Sciences,
Beijing Normal University, Beijing 100875, China}\\
\footnotesize{and}\\ \footnotesize{Department of Mathematics,
Swansea University, Singleton Park, SA2 8PP, UK}\\
\footnotesize{Email: wangfy@bnu.edu.cn; F.Y.Wang@swansea.ac.uk}}
\begin{document}
\def\R{\mathbb R}  \def\ff{\frac} \def\ss{\sqrt} \def\B{\mathbf
B}
\def\N{\mathbb N} \def\kk{\kappa} \def\m{{\bf m}}
\def\dd{\delta} \def\DD{\Delta} \def\vv{\varepsilon} \def\rr{\rho}
\def\<{\langle} \def\>{\rangle} \def\GG{\Gamma} \def\gg{\gamma}
  \def\nn{\nabla} \def\pp{\partial} \def\EE{\scr E}
\def\d{\text{\rm{d}}} \def\bb{\beta} \def\aa{\alpha} \def\D{\scr D}
  \def\si{\sigma} \def\ess{\text{\rm{ess}}}
\def\beg{\begin} \def\beq{\begin{equation}}  \def\F{\scr F}
\def\Ric{\text{\rm{Ric}}} \def\Hess{\text{\rm{Hess}}}
\def\e{\text{\rm{e}}} \def\ua{\underline a} \def\OO{\Omega}  \def\oo{\omega}
 \def\tt{\tilde} \def\Ric{\text{\rm{Ric}}}
\def\cut{\text{\rm{cut}}} \def\P{\mathbb P} \def\ifn{I_n(f^{\bigotimes n})}
\def\C{\scr C}      \def\aaa{\mathbf{r}}     \def\r{r}
\def\gap{\text{\rm{gap}}} \def\prr{\pi_{{\bf m},\varrho}}  \def\r{\mathbf r}
\def\Z{\mathbb Z} \def\vrr{\varrho} \def\ll{\lambda}
\def\L{\scr L}\def\Tt{\tt} \def\TT{\tt}\def\II{\mathbb I}
 \def\Sect{{\rm Sect}}\def\E{\mathbb E} \def\H{\mathbb H}
\def\M{\scr M}\def\Q{\mathbb Q} \def\texto{\text{o}} \def\LL{\Lambda}
\def\Rank{{\rm Rank}} \def\B{\scr B}   \def\HR{\hat{\R}^d}
\def\BB{\mathbb B}\def\vp{\varphi}

\maketitle
\begin{abstract} By solving a control problem and using Malliavin calculus, explicit derivative formula is derived for the semigroup $P_t$ generated by  the Gruschin type operator  on $\R^{m}\times \R^{d}:$
$$L (x,y)=\ff 1 2 \bigg\{\sum_{i=1}^m \pp_{x_i}^2 +\sum_{j,k=1}^d (\si(x)\si(x)^*)_{jk} \pp_{y_j}\pp_{y_k}\bigg\},\ \ (x,y)\in \R^m\times\R^d,$$ where $\si\in C^1(\R^m; \R^d\otimes\R^d)$ might be degenerate.   In particular,   if  $\si(x)$ is comparable with $|x|^{l}I_{d\times d}$ for some $l\ge 1$ in the sense of (\ref{A4}),  then for any $p>1$ there exists a constant $C_p>0$ such that $$|\nn P_t f(x,y)|\le \ff{C_p (P_t |f|^p)^{1/p}(x,y)}{\ss{t}\land \ss{t(|x|^2+t)^l}},\ \ t>0, f\in \B_b(\R^{m+d}), (x,y)\in \R^{m+d},$$ which implies  a new Harnack type inequality for the semigroup. A more general model is also investigated.
\end{abstract} \noindent

 AMS subject Classification:\ 60J75, 60J45.   \\
\noindent
 Keywords: Gruschin semigroup, derivative formula, gradient estimate.
 \vskip 2cm

\section{Introduction}

It is well-known that a hypoelliptic diffusion semigroup on $\R^d$ has a smooth transition density w.r.t. the Lebesgue measure (cf. \cite{N}). An interesting research topic is then to derive explicit estimates on the derivatives of the diffusion semigroup. To this end,  the derivative formula, which is called the Bismut formula or the Bismut-Elworthy-Li formula   due to \cite{B, EL}, has become a powerful tool.

 In the elliptic setting, the formula can be explicitly established by using the associated Bakry-Emery curvature tensor. But in the degenerate case the curvature is no-longer available and the existing formula established using the Malliavin covariance matrix is
  normally less explicit, see e.g.  \cite[Theorem 10]{AT0} and
 \cite[Theorem 3.2]{AT}. To establish  explicit derivative formulae for hypoelliptic semigroups, one has to build and solve  some  control problems   associated to the corresponding stochastic differential equations, see e.g. \cite{GW, WZ, Z} for the study of generalized stochastic Hamiltonian systems, and  see \cite[Section 6]{AT} for some simple examples. See also \cite{P} for the study of hypoelliptic Ornstein-Uhlenbeck semigroups.

Among Laplacian type hypoelliptic operators without drift term, two typical models are the Kohn-Laplacian on   Heisenberg groups and the Gruschin operator on $\R^2$. In recent years, the gradient estimate and applications have been intensively investigated for the heat semigroup $P_t$ generated by the Kohn-Laplacian on  finite- or infinite-dimensional Heisenberg groups, see \cite{BBBC, BGM, DM, Li} and the references within. In particular, the gradient inequality
\beq\label{1.1} \GG_1(P_t f)\le CP_t\GG_1(f),\ \ t\ge 0, f\in C_b^1, t\ge 0\end{equation} is confirmed in \cite{DM} for some constant $C>0$, where $\GG_1$ is the associated square field. This gradient inequality has important  applications, for instance, it implies the heat kernel Poincar\'e inequality and thus (cf. \cite{BBBC}),
\beq\label{1.1'} \GG_1(P_t f)\le \ff c t P_t  f^2,\ \ f\in C_b^1, t>0\end{equation}for some constant $c>0.$

 Accordingly,
 one may wish to prove (\ref{1.1}) and (\ref{1.1'}) also for the semigroup    generated by the  Gruschin operator $  \pp_x^2 +x^{2l} \pp_y^2$
 on $\R^2$, where $l\in \mathbb N.$  As pointed out to the author by the referee that when $l=1$ these can be confirmed by using the known inequalities on the Heisenberg group
 and the submersion $\psi: (x,y,z)\mapsto (x, z+\ff{xy}2).$ Indeed, letting $\tt P_t$ and $\tt \GG_1$ be the semigroup and square field associated to the Kohn-Laplacian $\tt X^2+\tt Y^2$ on $\R^3$, where $\tt X:= \pp_x-\ff y 2 \pp_x, \tt Y:= \pp_y +\ff x 2\pp_z$, we have
 $$  \tt X(f\circ \psi)= (\pp_x f)\circ\psi,\ \ \tt Y(f\circ \psi)= (x\pp_y f)\circ\psi,\ \ f\in C^1(\R^2),$$ so that
 $$\GG_1(P_t f)\circ\psi= \tt \GG_1(\tt P_t f\circ\psi),\ \ (P_t\GG_1(f))\circ\psi= \tt P_t\tt \GG_1(f\circ\psi)$$ hold.
 When $l\ge 2$, (\ref{1.1}) is however not yet available. We also would like to mention that for $l=1$,     the generalized curvature-dimension condition introduced and applied in \cite{BB1,BB2,BB3} holds, so that the corresponding  results, in particular the  gradient estimates and applications derived in \cite{BB2}, are valid. Even when $l\ge 2$,   although their generalized curvature condition is no longer available, a more general version of curvature condition has been confirmed in \cite{W12},  so that the $L^2$-gradient estimate as in Corollary \ref{C1.2} below for $p=2$ holds.

 In this paper, we aim to establish the Bismut-type derivative formula and gradient estimates for the semigroup generated by  the  following Gruschin-type operators on   $\R^{m+d}$:
$$L (x,y)=\ff 1 2 \bigg\{\sum_{i=1}^m \pp_{x_i}^2 +\sum_{j,k=1}^d (\si(x)\si(x)^*)_{jk} \pp_{y_j}\pp_{y_k}\bigg\},\ \ (x,y)\in \R^m\times\R^d=\R^{m+d},$$ where $\si\in C^1(\R^m; \R^d\otimes\R^d)$ might be degenerate.   In this general case, it seems hard to adopt the above mentioned arguments developed for Heisenberg groups and subelliptic operators satisfying the generalized curvature. Our study is based on     Malliavin calculus.

Let $\GG_1$ be the square field associated to $L$. Then
\beq\label{SQ} \GG_1(f)(x,y)= |\nn f(\cdot,y)(x)|^2 +|\si(x)^* \nn f(x,\cdot)(y)|^2,\ \ (x,y)\in \R^{m+d}, f\in C^1(\R^{m+d}).\end{equation}
  We will use $|\cdot|$ and $\|\cdot\|$ to denote the Euclidean norm and the operator norm   respectively.

To construct the associated diffusion process, we consider the  stochastic differential equation on $\R^{m+d}$:
\beq\label{E1} \beg{cases} \d X_t= \d B_t,\\
\d Y_t= \si(X_t) \d \tt B_t,\end{cases}\end{equation} where $(B_t,\tt B_t)$ is a Brownian motion on $\R^{m+d}$. It is easy to see that for any initial data the equation has a unique  solution and the solution is non-explosive. Let $\E^{x,y}$ stands for the expectation taken for the solution starting at $(x,y)\in\R^{m+d}$. We have
$$P_t f(x,y)= \E^{x,y} f(X_t, Y_t),\ \ \  f\in \B_b(\R^{m+d}), (x,y)\in \R^{m+d}, t\ge 0.$$
 To establish explicit derivative formula for $P_t$, we need the following assumption.
 \paragraph{(A)} For any $T>0$ and $x\in\R^m$, $Q_T:=  \int_0^T \si(x+B_t)\si(x+B_t)^*\d t$ is invertible such that
 $$\E\bigg\{\|Q_T^{-1}\|^2\int_0^T\big(\|\nn \si(x+B_t)\|^4+ \|\si(x+B_t)\|^4+1\big)\d t\bigg\}<\infty.$$
Obviously, $Q_T$ is invertible if so is $\si(x)$ for a.e. $x\in\R^m$. According to the proof of Corollary \ref{C1.2} below, assumption {\bf (A)}is ensured by (\ref{A4}) below.

\beg{thm}\label{T1.1} Assume {\bf (A)}.  For any $f\in C_b^1(\R^{m+d})$ and $v=(v_1,v_2)\in \R^{m+d}$,
$$\nn_v P_Tf(x,y)= \E^{x,y}\big\{f(X_T, Y_T) M_T\big\},\ \ (x,y)\in\R^{m+d}, T>0$$ holds for
\beg{equation*}\beg{split} M_T= &\ff {\<v_1, B_T\>}T  -  {\rm Tr}\bigg(Q_T^{-1} \int_0^T \ff{T-t}T\big\{(\nn_{v_1} \si)\si^*\big\} (x+B_t)\d t \bigg)\\
&+\bigg\<Q_T^{-1}\bigg\{v_2+ \int_0^T \ff{T-t}T(\nn_{v_1}\si)(x+B_t)   \d\tt B_t\bigg\}, \int_0^T \si(x+B_t) \d\tt B_t\bigg\>,\end{split}\end{equation*}where $\nn_v$ stands for the directional derivative along $v$.\end{thm}

To derive explicit estimates, we assume that $\si(x)$ is comparable with $|x|^lI_{d\times d}$ in the sense of (\ref{A4}) below.

\beg{cor} \label{C1.2} Let $l\in [1,\infty)$ and assume that
\beq\label{A4} \|\si(x)\|\ge a |x|^l,\ \ \|\si(x)\|+\|\nn \si(x)\|\cdot|x|\le b|x|^l,\ \ x\in\R^m\end{equation}holds for some constants $a,b>0$.
 Then for any $p>1$ there exists a constant $C_p>0$ such that for any $v=(v_1,v_2)\in\R^{m+d},$
\beq\label{A5}|\nn_v P_Tf(x,y)| \le  C_p(P_T|f|^p)^{1/p}(x,y)\bigg(  \ff {|v_1|} {\ss T}+\ff{|v_2|} {\ss{T(|x|^2+T)^l}}\bigg),\ \  T>0, (x,y)\in \R^{m+d}.\end{equation} Consequently,
\beq\label{A6}\GG_1(P_Tf)\le \ff{CP_Tf^2}{T},\ \ T>0, f\in \B_b(\R^{m+d})\end{equation} holds for some constant $C>0,$  where $\GG_1$ is given by $(\ref{SQ})$. \end{cor}

Let $P_t(z;\cdot)$ be the transition  probability kernel of $P_t$. It is easy to see that (\ref{A5}) implies
$$\|P_t((x,y); \cdot)-P_t((x',y'); \cdot)\|_{var} \le C\|f\|_\infty \ff{|x-x'|}{\ss T} + \ff{|y-y'|}{\ss{T^{3+l}}},\ \ T>0, (x,y), (x'y')\in\R^{m+d}$$ for some constant $C>0$,
 where $\|\varphi\|_{var}:= \sup \varphi(\cdot)-\inf \varphi(\cdot)$ is the total variational norm of a signed measure $\varphi$. Consequently (cf. \cite{LL}), the Markov process has successful couplings. Moreover, according to the following result, (\ref{A5}) and (\ref{A6}) also imply Harnack type inequalities for $P_T$.

In general, let $E$ be a connected differential manifold and let $\GG_1$ be a square field of type
$$\GG_1(f)=\sum_{i=1}^l (X_i f)^2$$ for some continuous vector  fields $\{X_i\}_{i=1}^d$. For any vector $v\in T_x E$, the intrinsic norm of $v$ induced by $\GG_1$ is
$$|v|_{\GG_1}=\sup\big\{|vf|(x): \GG_1(f)(x)\le 1\big\}.$$ For any $C^1$-curve $\gg: [0,1]\to E$, the length of $\gg$ induced by $\GG_1$ is
$$\ell(\gg)= \int_0^1 |\dot\gg_s|_{\GG_1}\d s.$$ Finally, for any $z,z'\in E$, the intrinsic distance   between them induced by $\GG_1$ is
$$\rr(z,z')= \inf\big\{\ell(\gg):\ \gg \text{\ is\ a\ } C^1\text{-curve\ linking\ } z\ \text{and}\ z'\big\}.$$ It is well known that $\rr$ is finite
if $\{X_i\}_{i=1}^d$ are smooth vector fields satisfying   H\"ormander's condition. An alternative way to define $\rr$ is to use the subunit curve. Recall that a $C^1$-curve $\gg: [0,T]\to E$ is called subunit w.r.t. $\GG_1$ if $|\ff{\d}{\d t} f(\gg_t)|\le \ss{\GG_1(f)(\gg_t)},\ t\in [0,T].$ Then
$$\rr(z,z')= \inf\big\{T>0:\ \text{there\ exists\ a\ subunit\ curve}\ \gg: [0,T]\to M, \gg_0=z, \gg_T=z'\big\}.$$

\beg{prp}\label{P1.3} Let $\GG_1$ and $\rr$ be fixed as above on a connected differential manifold $E$ such that $\rr$ is finite. Let $P$ be a $($sub-$)$Markov operator on $\scr B_b(E)$, the set of all bounded measurable functions on $E$. Then for any constant $C>0$,
\beq\label{A7} \GG_1(Pf)\le C^2Pf^2,\ \ f\in C_b^1(E)\end{equation} is equivalent to the Harnack type inequality
\beq\label{A8} Pf(z')\le Pf(z) +  C\rr(z,z')\ss{P f^2(z')},\ \ z,z'\in E, f\ge 0, f\in\B_b(E).\end{equation}  \end{prp}

A simple application of (\ref{A8}) is the following Harnack inequality for the transition kernel $P(z,\cdot)$ of $P$: taking $f=1_A$ in (\ref{A8})
for measurable set $A$, we obtain
$$P(z, \cdot)\le P(z',\cdot) +  C\rr(z,z')\ss{P(z',\cdot)},\ \ z,z'\in E.$$

We will prove Theorem \ref{T1.1}  in Section 2 and prove Corollary \ref{C1.2} and Proposition \ref{P1.3} in Section 3.   Finally, in section 4 we extend Theorem \ref{T1.1} to a more general model.

\section{Proof of Theorem \ref{T1.1} }

To establish the derivative formula,   we  first briefly recall the integration by parts formula for the Brownian motion. Let $T>0$ be fixed and let
$$\H=\bigg\{h\in C([0,T];\R^{m+d}):\ h(0)=0, \|h\|_\H^2:=\int_0^T  |h'(t)|^2 \d t<\infty\bigg\}$$ be the Cameron-Martin space.
Let $\mu$ be the distribution of $(B_t,\tt B_t)_{t\in [0,T]}$, which is a probability measure (i.e. Wiener measure) on the path space  $W=C([0,T];\R^{m+d})$.
A function $F\in L^2(W;\mu)$ is called differentiable if for any $h\in \H$, the directional derivative
$$D_h F:=\lim_{\vv\to 0} \ff{F(\cdot+\vv h)-F(\cdot)}\vv$$ exists in $L^2(W;\mu).$ We write $F\in \D(D)$ if moreover
$$\H\ni h \mapsto D_h F\in L^2(W;\mu)$$ is a bounded linear operator. In this case the Malliavin gradient $DF$ is defined as the unique element in $L^2(W\to\H;\mu)$  such that
$\<DF,h\>_\H=D_h F$ for $h\in\H$.  It is well known that $(D,\D(D))$ is a closed operator in $L^2(W;\mu)$, whose adjoint operator $(\dd,\D(\dd))$ is called the divergence
operator. That is,
\beq\label{INT} \int_W D_h F\d\mu= \int_W F \dd(h)\d\mu,\ \ \ F\in \D(D), h\in \D(\dd).\end{equation}

\beg{thm}\label{T2.1} For fixed $T>0$ and $v=(v_1,v_2)\in \R^{m+d}$, let $h_1\in C^1([0,T];\R^m)$ with $h_1(0)=0$ and $h_1(T)=v_1$. If there exists a process $\{h_2(t)\}_{t\in [0,T]}$ on $\R^d$ such that   $h_2(0)=0$,  and $h:=(h_1,h_2)\in \D(\dd)$ satisfying
\beq\label{2.1} \int_0^T \si(X_t) h_2'(t)\d t + \int_0^T (\nn_{h_1(t)-v_1}\si)(X_t)   \d\tt B_t =v_2,\end{equation} then
$$\nn_v P_T f= \E \big\{f(X_T,Y_T) \dd(h)\big\},\ \ f\in C_b^1(\R^2).$$\end{thm}

\beg{proof} From (\ref{E1}) it is easy to see that the derivative process $(\nn_v X_t, \nn_v Y_t)_{t\ge 0}$ solve the equation
$$\beg{cases} \d \nn_v X_t=0,\ \  & \nn_v X_0= v_1,\\
\d \nn_v Y_t = (\nn_{\nn_v X_t} \si)(X_t) \d\tt B_t,\ \ &\nn_v Y_0= v_2.\end{cases}$$ So,
\beq\label{nn} \beg{cases} \nn_v X_t= v_1,\\
\nn_v Y_t = v_2 +  \int_0^t(\nn_{v_1} \si)(X_s)  \d\tt B_s.\end{cases}\end{equation} Next, for $h$ given in the theorem, we have
$$\beg{cases} \d D_h X_t = h_1'(t) \d t,\ \ & D_h X_0=0,\\
\d D_h Y_t = \si(X_t) h_2'(t) \d t +(\nn_{D_h X_t}\si)(X_t)  \d\tt B_t,\ \ &D_h Y_0=0.\end{cases}$$ Thus,
$$\beg{cases} D_h X_t= h_1(t),\\
D_h Y_t = \int_0^t \si(X_s) h'_2(s) \d s + \int_0^t(\nn_{h_1(s)}\si)(X_s)  \d\tt B_s.\end{cases}$$ Since $h_1(T)=v_1$, combining this with   (\ref{2.1}) and (\ref{nn}) we obtain $$(\nn_v X_T, \nn_v Y_T)= (D_h X_T, D_h Y_T).$$  Therefore, for any $f\in C_b^1(\R^2)$, it follows from (\ref{INT}) that
\beg{equation*}\beg{split} \nn_v P_T f&= \E \<\nn f(X_T, Y_T), (\nn_v X_T, \nn_v Y_T)\>
 = \E  \<\nn f(X_T, Y_T), (D_h X_T, D_h Y_T)\>\\
 &=\E D_h \{f(X_T,Y_T)\}
=\E\{f(X_T,Y_T) \dd(h)\}.\end{split}\end{equation*}
\end{proof}

To prove Theorem \ref{T1.1}, the   key point is to solve the control problem (\ref{2.1}). To this end, we will   need the following fundamental lemma.

\beg{lem}\label{LL} Let   $\rr_t$ be a predictable process on $\R^d$ with $\E\int_0^T|\rr_t|^q<\infty$ for some $q\ge 2.$ Then
\beg{equation*}\beg{split} \E\bigg|\int_0^T \<\rr_t, \d\tt B_t\>\bigg|^q&\le \Big\{\ff{q(q-1)}2\Big\}^{q/2}\bigg(\int_0^T (\E|\rr_t|^q)^{2/q}\d t\bigg)^{q/2}\\
&\le \Big\{\ff{q(q-1)}2\Big\}^{q/2}T^{(q-2)/2} \int_0^T \E|\rr_t|^q\d t.\end{split}\end{equation*}\end{lem}

\beg{proof} It suffices to prove the first inequality since the second follows immediately from Jensen's inequality.
Let $N_t=\int_0^t \<\rr_s,\d\tt B_s\>,\ t\ge 0.$ Then $\d \<N\>_t= |\rr_t|^2\d t$ and
$$\d N_t^2 =2N_t\d N_t+ |\rr_t|^2\d t.$$ Noting that $|N_t|^q= (N_t^2)^{q/2}$, by It\^o's formula we obtain
\beg{equation*}\beg{split} \d |N_t|^q &= \ff q 2 (N_t^2)^{(q-2)/2}\d N_t^2 + \ff{q(q-2)}2 (N_t^2)^{(q-4)/2} N_t^2|\rr_t|^2\d t\\
&=qN_t|N_t|^{q-2}\d N_t +\ff{q(q-1)} 2 |N_t|^{q-2}|\rr_t|^2\d t.\end{split}\end{equation*} Therefore,
\beg{equation*}\beg{split} \E|N_T|^q &= \ff{q(q-1)} 2 \int_0^T \E\big\{|N_t|^{q-2} |\rr_t|^2\big\}\d t \\
&\le \ff{q(q-1)}2 \int_0^T \big(\E|N_t|^q\big)^{(q-2)/q} \big(\E|\rr_t|^q\big)^{2/q}\d t \\
&\le \ff{q(q-1)}2\big(\E|N_T|^q\big)^{(q-2)/q}  \int_0^T   \big(\E|\rr_t|^q\big)^{2/q}\d t.\end{split}\end{equation*} Up to an approximation argument we may assume that $\E|N_T|^q<\infty$, so that this implies
$$\E|N_T|^q \le \Big\{\ff{q(q-1)}2\Big\}^{q/2}\bigg(\int_0^T (\E|\rr_t|^q)^{2/q}\d t\bigg)^{q/2}.$$
\end{proof}

\beg{proof}[Proof of Theorem \ref{T1.1}] We  assume that $(X_0,Y_0)=(x,y)$ and
simply denote $\E^{x,y}$ by $\E$. Let
\beq\label{H1} h_1(t)= \ff {tv_1}T,\ \ t\in [0,T] \end{equation} and
\beq\label{H2} h_2(t)= \bigg(\int_0^t \si(X_s)^*\d s\bigg) Q_T^{-1} \bigg(v_2+  \int_0^T\ff{T-s}T(\nn_{v_1} \si)(X_s)  \d \tt B_s\bigg),\ \ t\in [0,T].\end{equation}
Then it is easy to see  that (\ref{2.1}) holds.  To see that $h:= (h_1, h_2)\in \D(\dd)$ and to calculate $\dd(h)$, let
\beg{equation*}\beg{split} &g_i= \bigg\<e_i, Q_T^{-1} \bigg(v_2+  \int_0^T\ff{T-s}T(\nn_{v_1} \si)(X_s) \d \tt B_s\bigg)\bigg\>,\\
&\tt h_i(t) =  \int_0^t \si(X_s)^*e_i \d s,\ \ \ i=1,\cdots, d,\end{split}\end{equation*} where $\{e_i\}_{i=1}^d$ is the canonical ONB on $\R^d$. We have
\beq\label{HH0}h(t)= (h_1(t), 0)+ \sum_{i=1}^d g_i  (0, \tt h_i(t)).\end{equation} It is easy to see that $h_1$ and $\tt h_i$ are adapted and
\beq\label{HH1} \beg{split} &\dd((h_1,0))= \int_0^T \<h_1'(t),\d B_t\>= \ff{\<v_1, B_T\>}T,\\
&\dd((0, \tt h_i))= \int_0^T \<\tt h_i'(t), \d\tt B_t\> = \int_0^T \<\si(X_t)^* e_i, \d\tt B_t\>. \end{split}\end{equation}
Let $\scr \C$ be the $\si$-field induced by $\{B_s: s\in [0,T]\}$. By Lemma \ref{LL} and noting that $X_t$ is measurable w.r.t. $\scr C$ while $\tt B$ is independent of $\scr C$, we have
\beg{equation*}\beg{split} &\E\Big(\big\{g_i \dd((0,\tt h_i))\big\}^2\Big|\scr C\Big) =\E\bigg(\bigg\{g_i\int_0^T \<\si(X_t)^* e_i, \d\tt B_t\>\bigg\}^2\bigg|\scr C\bigg)\\
&\le  2 \|v_2\|^2\|Q_T^{-1}\|^2 \E\bigg(\bigg\{\int_0^T \<\si(X_t)^*e_i, \d\tt B_t\>\bigg\}^2\bigg|\scr C\bigg)\\
 &\quad + 2 \E\bigg(\bigg\{\int_0^T \ff{T-t}T\big\<\big(\nn_{v_1}\si(X_t)\big)^*(Q_T^{-1})^* e_i, \d\tt B_t\big\>\bigg\}^2\bigg\{\int_0^T \<\si(X_t)^*e_i, \d\tt B_t\>\bigg\}^2\bigg|\scr C\bigg)\\
 &\le c \|Q_T^{-1}\|^2\int_0^T\|\si(X_t)\|^2\d t \\
 &\quad+ 2\bigg[\E\bigg(\bigg\{\int_0^T \ff{T-t}T\big\<\big(\nn_{v_1}\si(X_t)\big)^*(Q_T^{-1})^* e_i, \d\tt B_t\big\>\bigg\}^4
 \bigg|\scr C\bigg)\cdot\E\bigg(\bigg\{\int_0^T \<\si(X_t)^*e_i, \d\tt B_t\>\bigg\}^4\bigg|\scr C\bigg)\bigg]^{1/2}\\
 &\le c'   \|Q_T^{-1}\|^2   \int_0^T( \|\nn \si(X_t)\|^4+ \|\si(X_t)\|^4+1)\d t \end{split}\end{equation*}  for some constants $c,c'>0$. So,   {\bf (A)} implies   $g_i \dd((0,\tt h_i))\in L^2(\P)$ for $i=1,\cdots d.$ Hence,   if for any $i\in\{1,\cdots d\}$ one has   $D_{(0,\tt h_i)}g_i\in L^2(\P)$, then $h\in\D(\dd)$ and by (\ref{HH0}) and (\ref{HH1}),
\beq\label{*W} \dd(h) = \ff{\<v_1, B_T\>}T +\sum_{i=1}^d\bigg\{g_i \int_0^T \<\si(X_t)^*e_i,\d\tt B_t\>-D_{(0,\tt h_i)}g_i\bigg\}.\end{equation}
Noting that $X_t=x+B_t$ is independent of $\tt B$, it is easy to see that
\beg{equation*}\beg{split} D_{(0,\tt h_i)} g_i &= \bigg\<e_i, Q_T^{-1} \int_0^T \ff{T-t}T (\nn_{v_1}\si)(X_t) \tt h_i'(t)\d t\bigg\>\\
&= \bigg\<e_i, Q_T^{-1} \int_0^T \ff{T-t}T \big\{(\nn_{v_1}\si)  \si^*\big\}(X_t)e_i\d t\bigg\>,\end{split}\end{equation*} which is in $L^2(\P)$ according to {\bf (A)}. Combining this with (\ref{*W}) and noting that $X_t=x+B_t$, we conclude that $h\in \D(\dd)$ and
$\dd(h)=M_T$.     Then the proof is finished by  Theorem \ref{T2.1}.
\end{proof}

\section{Proofs of Corollary \ref{C1.2} and Proposition \ref{P1.3}}
To verify {\bf(A)} for $\si$  given in Corollary \ref{C1.2}, we first present the following lemma.
\beg{lem} \label{L3.1}   For any $n\in [1,\infty)$ and $\aa>0$, there exists a constant $c>0$ such that
$$\E^{x,y}\bigg(\int_0^T |X_t|^{2n}\d t\bigg)^{-\aa}\le \ff c {T^\aa (|x|^2+T)^{\aa n}},\ \ \ T>0, (x,y)\in\R^2. $$ \end{lem}

\beg{proof} We shall simply denote $\E^{x,y}$ by $\E$. Since $X_t= x+B_t$, for any $\ll>0$ we have (see e.g. \cite[page 142]{Hand})
 \beg{equation*} \beg{split}&\E\e^{-\ll \int_0^T  |X_t|^2\d t}=\prod_{i=1}^m \E \e^{-\ll \int_0^T (x_i+B_t^{(i)})^2\d t} \le \ff {\exp[-\ff{x^2\ss\ll}{\ss 2}\tanh(\ss{2\ll}\,T)]} {\{\coth\big(\ss {2\ll} T\big)\}^{m/2}}\\
 & \le 2^{m/2}\, \exp\bigg[ -\ff{mT\ss\ll}{\ss 2} - \ff{x^2\ss\ll}{2\ss 2} \Big\{\big(\ss{2\ll}\, T\big)\land 1\Big\}\bigg]\\
 &\le 2^{m/2} \, \exp\bigg[  - \ff{(x^2+T)\ss\ll}{2\ss 2} \Big\{\big(\ss{2\ll}\, T\big)\land 1\Big\}\bigg]\\
 &\le 2^{m/2} \exp\bigg[-\ff{(x^2+T)\ss\ll}{2\ss 2}\bigg]+2^{m/2} \exp\bigg[-\ff{(x^2+T)\ll T}{\ss 2}\bigg].\end{split}\end{equation*}
 This implies that for any $r>0,$

\beg{equation*}\beg{split}  & \E\exp\bigg[-\ll\int_0^T |X_t|^{2n}\d t\bigg]= \E \exp\bigg[- \int_0^T \big(\ll^{1/n} |X_t|^2\big)^n  \d t\bigg]   \\
&\le    \E \exp\bigg[\int_0^T \Big(\ff{n-1}{n^{n/(n-1)}}r^{n/(n-1)} -r \ll^{1/n} |X_t|^2\Big)\d t\bigg]\\
&\le  2^{m/2}\, \exp\bigg[ \ff{T(n-1)}{n^{n/(n-1)}} r^{n/(n-1)}\bigg]\bigg(\exp\bigg[   -\ff{(x^2+T) \ll^{1/(2n)}\ss r} {2\ss 2} \bigg]+\exp\bigg[ -\ff{(x^2+T)T \ll^{1/n}r} {\ss 2}   \bigg]\bigg).\end{split}\end{equation*} Taking $r= T^{-(n-1)/n}$ we obtain
$$\E\exp\bigg[-\ll\int_0^T |X_t|^{2n}\d t\bigg]\le c_1\bigg(\exp\bigg[-\ff{(x^2+T)\ll^{1/(2n)}}{2\ss 2\, T^{(n-1)/2n}}\bigg]+ \exp\bigg[-\ff{(x^2+T)(\ll T)^{1/n}}{\ss 2}\bigg]\bigg)$$ for some constant $c_1>0$.  Noting that
 $$\int_0^\infty \ll^{\aa-1}\e^{-\theta \ll^{1/l}}\d\ll= \ff{l}{\theta^{\aa l}}\int_0^\infty \e^{-s}s^{\aa l-1} \d s=\ff{l \GG(\aa l)}{\theta^{\aa l}}$$ holds for all
 $l\ge 1$ and $\theta,\aa>0$, we conclude that
\beg{equation*}\beg{split} &\E\bigg(\int_0^T|X_t|^{2n}\d t\bigg)^{-\aa} =\ff 1 {\GG(\aa)} \int_0^\infty \ll^{\aa-1} \E \exp\bigg[-\ll\int_0^T |X_t|^{2n}\d t\bigg]\d\ll\\
&\le c_1\int_0^\infty \ll^{\aa-1}  \bigg\{\exp\bigg[-\ff{(|x|^2+T)\ll^{1/(2n)}}{2\ss 2\, T^{(n-1)/2n}}\bigg]+ \exp\bigg[-\ff{(|x|^2+T)(\ll T)^{1/n}}{\ss 2}\bigg]
\bigg\}\d\ll\\
&\le   \ff{c_2T^{\aa(n-1)}}{(|x|^2+T)^{2\aa n}}+\ff {c_3}{(|x|^2+T)^{\aa n}T^\aa}\le \ff{c}{(|x|^2+T)^{\aa n}T^\aa}\end{split}\end{equation*} holds for some constants $c_2,c_3$ and $c$.
\end{proof}

\beg{proof}[Proof of Corollary \ref{C1.2}]  By   Jensen's inequality, it suffices to prove for $p\in (1,2]$ so that
$q:= \ff p{p-1}\ge 2.$ In fact, once  (\ref{A5}) holds for $p=2$, it also holds for $p>2$ with $C_p=C_2$ since in this case
$(P_Tf^2)^{1/2}\le (P_T|f|^p)^{1/p}.$

It is easy to see that (\ref{A4}) implies
$$Q_T\ge \bigg(a^2 \int_0^T |X_t|^{2l}\d t\bigg)I_{d\times d},$$ and hence,
\beq\label{AB} \|Q_T^{-1}\|\le\ff 1 {a^2\int_0^T|X_t|^{2l}\d t}.\end{equation} Since $\{X_t\}_{t\in [0,T]}$ is measurable w.r.t. $\scr C$ and due to (\ref{A4})
$$\|\{(\nn_{v_1}\si)\si^*\}(X_t) \|\le b^2 |X_t|^{2l-1},$$  we obtain
\beq\label{GG1} \beg{split} &\E\bigg(\bigg|\ff{\<v_1,B_T\>}T- \text{Tr}\bigg(Q_T^{-1}\int_0^T \ff{T-t}T \big\{(\nn_{v_1}\si)\si^*\big\} (X_t)\d t\bigg)  \bigg|^q\bigg|\scr C\bigg)\\
&=\bigg|\ff{\<v_1,B_T\>}T- \text{Tr}\bigg(Q_T^{-1}\int_0^T \ff{T-t}T \big\{(\nn_{v_1}\si) \si^*\big\} (X_t)\d t\bigg)  \bigg|^q\\
&\le c_1 |v_1|^q \bigg(\ff{|B_T|^q}{T^q} + \ff{T^{q-1} \int_0^T |X_t|^{(2l-1)q}\d t}{(\int_0^T|X_t|^{2l}\d t)^q}\bigg)\end{split}\end{equation} for some constant $c_1>0.$
Moreover, since    $\tt B_t$ is independent of $\scr C$, due to (\ref{AB}) and Lemma \ref{LL} there exist constants $c_2,c_3>0$ such that
$$\E\bigg(\|Q_T^{-1}\|^q\bigg|\Big\<v_2, \int_0^T\si(X_t)\d\tt B_t\Big\>\bigg|^q\bigg|\scr C\bigg)\le \ff{c_2 |v_2|^qT^{q/2-1} \int_0^T |X_t|^{lq}\d t}{ (\int_0^T|X_t|^{2l}\d t)^q}$$ and
\beg{equation*}\beg{split} &\E\bigg(\|Q_T^{-1}\|^q\bigg|\Big\<\int_0^T\ff{T-t}T(\nn_{v_1}\si)(X_t)   \d\tt B_t,\ \int_0^T \si(X_t)\d\tt B_t\>\bigg|^q\bigg|\scr C\bigg)\\
&\le \ff{c_2}{(\int_0^T |X_t|^{2l}\d t)^q} \bigg\{\E\bigg(\bigg|\int_0^T\ff{T-t}T (\nn_{v_1}\si)(X_t) \d\tt B_t\bigg|^{(2l-1)q/(l-1)}\bigg|\scr C\bigg)\bigg\}^{(l-1)/(2l-1)}\\
&\qquad \times
\bigg\{\E\bigg(\bigg|\int_0^T\si(X_t)\d\tt B_t\bigg|^{(2l-1)q)/l}\bigg|\scr C\bigg)\bigg\}^{l/(2l-1)}\\
&\le \ff{c_3|v_1|^qT^{q-1}\int_0^T|X_t|^{(2l-1)q}\d t}{(\int_0^T |X_t|^{2l}\d t)^q}\end{split}\end{equation*} hold. Combining these with (\ref{GG1}) we obtain
\beg{equation}\label{GG2}\beg{split}& \E|M_T|^q= \E\big\{\E(|M_T|^q|\scr C)\big\}\\
&\le c_3  \E \bigg\{\ff {|v_1|^q} {T^{q/2}} + \ff{|v_1|^qT^{q-1}\int_0^T|X_t|^{(2l-1)q}\d t}{(\int_0^T|X_t|^{2l}\d t)^q} + \ff {|v_2|^qT^{\ff q 2 -1}\int_0^T|X_t|^{lq}\d t} { (\int_0^T|X_t|^{2l}\d t)^q}\bigg\}.\end{split}\end{equation} By Lemma \ref{L3.1} and noting that $X_t= x+B_t$, we conclude  that for any $\bb\ge 1$,
\beg{equation*}\beg{split} &\E\bigg\{\ff{\int_0^T|X_t|^\bb \d t }{(\int_0^T|X_t|^{2l}\d t)^q} \bigg\}\le \bigg\{\E\bigg(\int_0^T|X_t|^\bb\d t\bigg)^2\bigg\}^{1/2}
\bigg\{\E\bigg(\int_0^T|X_t|^{2l}\d t\bigg)^{-2q}\bigg\}^{1/2}\\
&\le \ff{c_3(T\E\int_0^T|X_t|^{2\bb}\d t)^{1/2}}{T^q(|x|^2+T)^{ql}}\le \ff{c_4T(|x|^2+T)^{\bb/2}}{T^q(|x|^2+T)^{ql}}=\ff{c_4}{T^{q-1}(|x|^2+T)^{ql-\bb/2}}\end{split}\end{equation*} holds for some constants $c_3,c_4>0.$ Substituting this into
(\ref{GG2})  we arrive at
$$ (\E|M_T|^q)^{1/q} \le c_5 \bigg\{\ff {|v_1|} {T^{q/2}} +  \ff {|v_2|} {T^{q/2}(|x|^2+T)^{ql/2}}\bigg\}^{1/q} $$  for some constant $c_5 >0.$ Therefore, (\ref{A5}) follows   since according to Theorem \ref{T1.1}
$$|\nn_v P_Tf(x,y)|= |\E\{f(X_T,Y_T)M_T\}|\le (P_T |f|^p)^{1/p} (\E|M_T|^q)^{1/q}.$$
\end{proof}

\beg{proof}[Proof of Proposition \ref{P1.3}]  {\bf (\ref{A7}) $\Rightarrow$ (\ref{A8})}. By the monotone class theorem, it suffices to prove (\ref{A8}) for $f\in C_b(E)$. For $z,z'\in E$, let $\rr=\rr(z,z').$ Up to an approximation argument we assume  that $\rr$ is reached by a subunit curve   $\gg: [0,\rr]\to E$ with $\gg_0=z,\gg_\rr=z'.$
 Then, due to (\ref{A7}),  for any positive $f\in C_b(E)$  we have
\beg{equation*}\beg{split} \ff{\d}{\d s} P\Big(\ff{f}{1+rs f}\Big)(\gg_s) &\le- P\Big(\ff{r f^2}{(1+rs f)^2}\Big)(\gg_s) +\ss{\GG_1\Big(\ff{f}{1+rs f}\Big)(\gg_s)}  \\
& \le -rP\Big(\ff{f^2}{(1+rs f)^2}\Big)(\gg_s) + C  \ss{P\Big(\ff f{1+rs f}\Big)^2(\gg_s)}\\
&\le \ff{C^2 }{4 r}.\end{split}\end{equation*}
Integrating over $[0,\rr]$ w.r.t. $\d s$ we obtain
$$P \Big(\ff f {1+r\rr f}\Big)(z')\le P f(z) +\ff{C^2\rr } {4r}.$$ Combining this with the fact that
$$\ff f{1+r\rr f} =f -\ff{r\rr f^2}{1+r\rr f} \ge f -r\rr f^2,$$ we obtain
$$Pf (z')\le Pf(z) +\ff{C^2 \rr  } {4r} + r\rr Pf^2(z').$$ Minimizing the right-hand side in $r>0$ we prove (\ref{A8}).

{\bf (\ref{A8}) $\Rightarrow$ (\ref{A7})}. By (\ref{A8}), we have
$$|Pf (x)-Pf(z')|\le  C\rr(x,y) \|f\|_\infty,\ \ f\in C_b(M).$$ So, $Pf$ is $\rr$-Lipschitz continuous for any $f\in \B_b(E)$. Let $z\in E$ and
$\gg: [0,1]\to M$ be $C^1$-curve such that $\gg_0=z, \rr(\gg_0,\gg_s)=s$ and $$\ff{\d}{\d s}Pf(\gg_s)|_{s=0} =\ss{\GG(Pf)(z)}.$$
Then it follows from (\ref{A8}) that
$$ \ss{\GG(Pf)(z)}= \lim_{s\to 0} \ff{P f(\gg_s)- Pf(\gg_0)} s \le C\lim_{s\to 0} \ss{Pf^2(\gg_s)} = C\ss{Pf^2(z)}.$$ Therefore, (\ref{A7}) holds.
\end{proof}

\section{An extension}

Consider the following SDE on $\R^{m+d}$:
\beq\label{E2} \beg{cases} \d X_t= \si_1(X_t)\d B_t+b_1(X_t)\d t,\\
\d Y_t= \si_2(X_t) \d \tt B_t+b_2(X_t)\d t,\end{cases}\end{equation} where $(B_t,\tt B_t)$ is a Brownian motion on $\R^{m+d}$, $\si_1\in C^1_b(\R^m; \R^m\otimes\R^m)$ is invertible with $\|\si_1^{-1}\|\le c$ for some constant $c>0$, $\si_2\in C^1(\R^m;\R^d\otimes\R^d)$ might be degenerate, $b_1\in C_b^1(\R^m;\R^m)$ and $b_2\in C^1(\R^m;\R^d).$  It is easy to see that for any initial data the solution exists uniquely  and is non-explosive. Let $P_t$ be the associated Markov semigroup. To establish the derivative formula, let $v=(v_1,v_2)\in \R^{m+d}$ and $T>0$ be fixed, and let $\xi_t$ solve the following SDE on $\R^m$:
\beq\label{E3} \d\xi_t =(\nn_{\xi_t}\si_1)(X_t)\d B_t +\Big\{(\nn_{\xi_t}b_1)(X_t)-\ff{\xi_t}{T-t}\Big\}\d t,\ \ \xi_0=v_1.\end{equation}
Since $\nn \si$ and $\nn b_1$ are bounded, the equation has a unique solution up to time $T$. It is easy to see from the It\^o formula   that
\beg{equation*}\beg{split} \d\bigg\{\ff{|\xi_t|^2}{T-t}\bigg\}&= 2\Big\<\ff{\xi_t}{T-t}, (\nn_{\xi_t}\si)(X_t)\d B_t\Big\> +\bigg(\ff{\|(\nn_{\xi_t}b_1)(X_t)\|^2}{T-t}+\ff{2\<\xi_t,(\nn_{\xi_t}b_1)(X_t)\>}{T-t}-\ff{|\xi_t|^2}{(T-t)^2}\bigg)\d t\\
&\le 2\Big\<\ff{\xi_t}{T-t}, (\nn_{\xi_t}\si)(X_t)\d B_t\Big\> +\Big(\ff{C|\xi_t|^2}{T-t}-\ff{|\xi_t|^2}{(T-t)^2}\Big)\d t,\ \ t\in [0,T) \end{split}\end{equation*} holds for some constant $C>0$. This implies that for $t\in [0,T)$,
\beq\label{E4} \E |\xi_t|^2 \le (T-t)\e^{Ct},\ \ \E\int_0^T \ff{|\xi_t|^2}{(T-t)^2}\d t<\infty,\ \ \end{equation}
Consequently, we may set $\xi_T=0$ so that $ \xi_t $ solves (\ref{E3}) for $t\in [0,T].$ Moreover,  for any $n\ge 1$ we have
$$\d  |\xi_t|^{2n} \le 2 n |\xi_t|^{2(n-1)}\Big\<\xi_t, (\nn_{\xi_t}\si)(X_t)\d B_t\Big\> +c(n)|\xi_t|^{2n}\d t$$ for some constant $c(n)\ge 0$. Therefore,
\beq\label{E5} \sup_{t\in [0,T]}\E |\xi_t|^{2n}<\infty,\ \ \ n\ge 1.\end{equation}
We are now able to state the derivative formula for $P_t$ as follows.

\beg{thm} Let $Q_T=\int_0^T\si_2(X_t)\si_2(X_t)^*\d t$ be invertible such that
\beq\label{E6} \E^{x,y} \bigg(\|Q_T^{-1}\|^2\int_0^T\big\{\|\si_2(X_t)\|^4+\|\nn\si_2(X_t)\|^4+\|\nn b_2(X_t)\|^4+1\big\}\d t\bigg)<\infty.\end{equation} Then
$$\nn_v P_T f(x,y)= \E^{x,y}\big\{f(X_T,Y_T)M_T\big\}$$ holds for $f\in C_b^1(\R^{m+d})$ and
\beg{equation*}\beg{split} M_T= &\int_0^T \Big\<\ff{\si_1(X_t)^{-1}\xi_t}{T-t}, \d B_t\Big\>-  {\rm Tr}\bigg(Q_T^{-1} \int_0^T \ff{T-t}T
\big\{(\nn_{\xi_t} \si_2)\si_2^* \big\}  (X_t)\d t \bigg)\\
&+\bigg\<Q_T^{-1}\bigg\{v_2+ \int_0^T \ff{T-t}T (\nn_{\xi_t}\si_2)(X_t)    \d\tt B_t+\int_0^T (\nn_{\xi_t} b_2)(X_t)\bigg\}, \int_0^T \si_2(X_t) \d\tt B_t\bigg\>.\end{split}\end{equation*}\end{thm}

\beg{proof} Let $h=(h_1,h_2)$, where
\beg{equation*}\beg{split} &h_1(t)= \int_0^t \ff{\si_1(X_s)^{-1}\xi_s}{T-s}\,\d s,\ t\in [0,T],\\
&h_2(t)= \bigg(\int_0^t \si_2(X_s)^*\d s\bigg)Q_T^{-1} \bigg(v_2 +\int_0^T( \nn_{\xi_t}\si_2)(X_t)\d\tt B_t+\int_0^T (\nn_{\xi_t}b_2)(X_t)\d t\bigg).\end{split}\end{equation*}
As in the proof of Theorem \ref{T1.1}, it is easy to see from (\ref{E4}), (\ref{E5}), (\ref{E6}) and $\|\si_1^{-1}\|\le c$ that $h\in \D(\dd)$ with
$\dd(h)=M_T.$ Therefore, it remains to verify that $(\nn_v X_T,\nn_v Y_T)=(D_h X_T, D_h Y_T).$ It is easy to see that both $\nn_v X_t$ and $D_h X_t+\xi_t$ solve the equation
$$\d V_t= (\nn_{V_t}\si_1)(X_t) \d B_t + (\nn_{V_t} b_1)(X_t)\d t,\ \ t\in [0,T], V_0=  v_1.$$ By the uniqueness of the solution we have $\nn_v X_t= D_h X_t+\xi_t$ for $t\in [0,T].$ Since $\xi_T=0,$ this implies that $\nn_v X_T= D_h X_T$. Moreover, we have
$$\beg{cases} \d\nn_v Y_t= (\nn_{\nn_v X_t}\si_2)(X_t)\d\tt B_t + (\nn_{\nn_v X_t}b_2)(X_t)\d t,\ \ \nn_vY_0=v_2,\\
\d D_h Y_t = (\nn_{D_h X_t}\si_2)(X_t)\d\tt B_t +\si_2(X_t)h_2'(t)\d t +(\nn_{D_h X_t}b_2)(X_t)\d t,\ \ D_h Y_0= 0.\end{cases}$$ Combining this with the definition of $h_2$ and $D_h X_t= \nn_v X_t-\xi_t$, we obtain
$$ D_h Y_T = \nn_v Y_T -v_2 -\int_0^T (\nn_{\xi_t}\si_2)(X_t)\d\tt B_t +\int_0^T \si_2(X_t) h_2'(t)\d t-\int_0^T (\nn_{\xi_t} b_2)(X_t)\d t= \nn_v Y_T.$$
Therefore, the proof is finished. \end{proof}

\paragraph{Acknowledgement.} The author would like to thank the referee for careful reading and   useful comments.

\beg{thebibliography}{99}

\bibitem{AT0}   M. Arnaudon, A. Thalmaier,
  \emph{Bismut type differentiation of semigroups,} Prob. Theory and Math. Stat. (Vilnius, 1998), 23-32,  VSP/TEV, Utrecht and Vilnius, 1999.

\bibitem{AT}   M. Arnaudon, A. Thalmaier,
  \emph{The differentiation of hypoelliptic diffusion semigroups,}    arXiv:1004.2174.

\bibitem{ATW09} M. Arnaudon, A. Thalmaier, F.-Y. Wang,
  \emph{Gradient estimates and Harnack inequalities on non-compact Riemannian manifolds,}
   Stoch. Proc. Appl. 119(2009), 3653--3670.

   \bibitem{BBBC} D. Bakry,  F. Baudoin,  M. Bonnefont, D. Chafa\"{\i}, \emph{ On gradient bounds for the heat kernel on the Heisenberg group,} J. Funct. Anal. 255 (2008),   1905--1938.

\bibitem{BB1} F. Baudoin, M. Bonnefont, \emph{Log-Sobolev inequalities for subelliptic operators satisfying a generalized curvature dimension inequality,} arXiv:1106.0491.12

\bibitem{BB2} F. Baudoin, M. Bonnefont, N. Garofalo, \emph{A sub-Riemannian Curvature-dimension inequality, volume doubling property and the Poincar\'e inequality,}
 arXiv:1007.1600.

\bibitem{BB3} F. Baudoin, N. Garofalo, \emph{Curvature-dimension inequalities and Ricci lower bounds for sub-Riemannian manifolds with transverse symmetries,} arXiv:1101.3590.

\bibitem{BGM} F. Baudoin, M. Gordina, T. Melcher, \emph{Quasi-invariance for heat kernel measures on sub-Riemannian infinite-dimensional Heisenberg groups,} arXiv:1108.1527.

\bibitem{B} J. M. Bismut, \emph{Large Deviations and the
Malliavin Calculus,} Boston: Birkh\"auser, MA, 1984.

\bibitem{Hand} A. N. Borodin, P. Salminen, \emph{Handbook of Brownian Motion - Facts and Formulae,} Birkh\"auser, Berlin, 1996.

\bibitem{EL}  K.D. Elworthy  and  X.-M. Li, \emph{Formulae for the
derivatives of heat semigroups,} J. Funct. Anal. 125(1994),
252--286.

\bibitem{DM} B. K. Driver, T. Melcker, \emph{Hypoelliptic heat kernel inequalities on the Heisenberg group,} J. Funct. Anal. 221(2005), 340--365.

\bibitem{GW} A. Guillin, F.-Y. Wang,
  \emph{Degenerate Fokker-Planck equations : Bismut formula, gradient estimate  and Harnack inequality,} to appear in J. Diff. Equat.
    arXiv:1103.2817

\bibitem{Li} H.-Q. Li, \emph{Estimation optimale du gradient du semi-groupe de la chaleur sur le groupe de Heisenberg,} J. Funct. Anal. 236(2006), 369--394.

\bibitem{LL} T. Lindvall, \emph{Lectures on the Coupling Methods,}  Wiley, New York, 1992.

\bibitem{N} D. Nualart, \emph{The Malliavin Calculus and Related Topics,} Springer, Berlin, 2006.

\bibitem{P}  E. Priola, \emph{ Formulae for the derivatives of
degenerate diffusion semigroups,}   J. Evol. Equ.  {\bf6} (2006),
557--600.

\bibitem{W12} F.-Y. Wang, \emph{Generalized curvature condition for subelliptic diffusion processes,} arXiv: 1202.0778.

\bibitem{WZ} F.-Y. Wang, X.-C. Zhang, \emph{Derivative formula and applications  for degenerate diffusion semigroups,} arXiv1107.0096.

\bibitem{Z} X.-C. Zhang, \emph{Stochastic flows and Bismut formulas for stochastic Hamiltonian systems,} Stoch. Proc. Appl. 120(2010), 1929--1949.

\end{thebibliography}
\end{document}